\documentclass[12pt]{article}
\usepackage{setspace}

\usepackage{amsfonts,pdfsync}
\usepackage{amsthm}
\usepackage{amsmath}
\usepackage{cases}

\usepackage{graphicx}
\usepackage{grffile}
\usepackage{adjustbox}

\usepackage{comment}
\usepackage[usenames,dvipsnames]{xcolor}
\usepackage{enumerate}
\usepackage{extarrows}

\usepackage{stackrel}
\usepackage{centernot}
\usepackage{caption}
\usepackage{subcaption}
\usepackage{fancyhdr}

\let\quoteOLD\quote
\def\quote{\quoteOLD\small}

\definecolor{labelkey}{cmyk}{0,0.8,1,0.5}
\definecolor{refkey}{cmyk}{0,0.8,1,0.5}

\newtheorem{theorem}{Theorem}[section]
\newtheorem{example0}{\sc Example}[subsection]

\newtheorem{proposition}{Proposition}

\newtheorem{lemma}{Lemma}

\numberwithin{equation}{section}
\numberwithin{theorem}{section}
\numberwithin{corollary}{section}
\numberwithin{proposition}{section}
\numberwithin{lemma}{section}
\numberwithin{definition}{section}
\numberwithin{remark}{section}

\makeatletter
\def\th@newremark{\th@remark\thm@headfont{\bfseries}}
\makeatletter

\def\boxit#1{\vbox{\hrule\hbox{\vrule\kern6pt
          \vbox{\kern6pt#1\kern6pt}\kern6pt\vrule}\hrule}}

\newcommand{\Fbar}{\overline{F}}
\newcommand{\Hbar}{\overline{H}}

\newcommand{\Gbar}{\overline{G}}

\newcommand{\topr}{\stackrel{\mathrm{P}}{\longrightarrow}}
\newcommand{\todr}{\stackrel{\mathrm{D}}{\longrightarrow}}

\newcommand{\R}{\Bbb{R}}

\newcommand{\N}{\Bbb{N}}

\newcommand{\rmd}{{\rm d}}

\newcommand{\wt}{\widetilde}

\newcommand{\dto}{\downarrow}

\newcommand{\be}{\begin{equation}}
\newcommand{\ee}{\end{equation}}
\newcommand{\bea}{\begin{eqnarray}}
\newcommand{\eea}{\end{eqnarray}}
\newcommand{\bean}{\begin{eqnarray*}}
\newcommand{\eean}{\end{eqnarray*}}
\newcommand{\ben}{\begin{equation*}}
\newcommand{\een}{\end{equation*}}
\newcommand{\ba}{\begin{align}}
\newcommand{\ea}{\end{align}}

\def\nexto{\kern -0.54em}

\usepackage[round, authoryear]{natbib}
\date{}

\begin{document}
\bibliographystyle{plainnat}

\title{Exact and Asymptotic Tests for Sufficient Followup in Censored Survival Data}


\author{Ross Maller$^a$,  Sidney Resnick$^b$ and Soudabeh Shemehsavar$^{a,b}$
\thanks{
$^a$Research School of Finance, Actuarial Studies \& Statistics, 
Australian National University, Ross.Maller@anu.edu.au;  
\newline
$^b$School of Operations Research and Information Engineering, Cornell University,
sir1@cornell.edu; 
\newline
$^c$School of Mathematics, Statistics \& Computer Sciences, University of Tehran, 
\newline
shemehsavar@ut.ac.ir (corresponding author)
}}


\maketitle

{}\null
\vskip-2.5cm 

\begin{onehalfspacing}

\begin{abstract}
\noindent
The existence of immune or cured individuals in a population and whether there is sufficient follow-up in a sample of censored observations on their lifetimes to be confident of their presence are questions of major importance in medical survival analysis.
So far only a few candidates  have been put forward as possible test  statistics for the existence of sufficient follow-up in  a sample. 
Here we investigate one such statistic, $Q_n$,  and give a detailed analysis
assuming  independence between survival and censoring times.
We obtain an exact finite sample as well as asymptotic distributions for $Q_n$,  and use these to calculate the power of the test as a function of the follow-up in the sample.
A particularly useful finding is that the asymptotic distribution of the  test statistic is parameter free in the null case when follow-up is insufficient.
The methods are illustrated with detailed schematic and real data sets, and the effect of dependence between survival and censoring is considered via a copula model and simulations.
\end{abstract}

\noindent{\bf Keywords} Sufficient follow-up;  censored survival  data;   cure model.

\end{onehalfspacing}

%
%

\section{Introduction}\label{intro}
There is a large and growing interest in the analysis of  censored survival data from a 
 population which may contain  immune or cured individuals, that is,  individuals who will not experience the event of interest no matter how long follow-up may be.
 A systematic formulation and treatment of this kind of problem is in  the book by
\cite{maller:zhoubook:1996} 
which contains  many practical examples from medicine, criminology and various other fields,  of this kind of data.
For other reviews and applications, see for example
 \cite{av:2018}, 
\cite{emvz:2020}, 
\cite{legrand:bertrand:2019},
 \cite{Ma:2009}, 
\cite{oblc:2012}, 
\cite{pt:2014},
\cite{pb:2021}.

 
In a sample of data of the kind mentioned, we have observations on the time to an event of interest (we refer to them as ``lifetimes"), possibly right-censored, of individuals from a population which contains some who are ``susceptible" to suffering the event under consideration, and possibly also some who are ``immune to", ``cured" of it, or are ``longterm survivors".
For ``susceptibles", the issue is to infer properties of their lifetime distribution from the sample. 
For ``immunes",  the main questions of interest are whether they are in fact present in the population, and if so in what proportion. Herein, we concentrate on the first question: based on the sample information, how confident can we be that immunes are in fact present in the population?
Currently developed methods for assessing this ultimately rely in some way on the amount of follow-up in the sample.

We do not know whether a particular censored lifetime in the sample is from an  immune or cured individual (uncensored lifetimes are obviously not from immunes); but, in aggregate,  the  presence of cured individuals may be signalled by an interval of constancy of the Kaplan-Meier estimator (\cite{kaplan:meier:1958}   (KME))  at its right hand end; 
 that is, the interval containing the censored lifetimes exceeding the largest uncensored lifetime. 
The length of that  interval and the  number of censored survival times larger than the largest uncensored survival time
  are important indicators for the presence of cured individuals, and  for whether there is sufficient follow-up in the sample to be confident of their presence.

Ways of testing for sufficient follow-up are still in  a very undeveloped state. 
One such test statistic,  $Q_n$, is suggested in \cite{MZ1994} and  \cite{maller:zhoubook:1996}, p.81.  
Since then, there have been only  two other definite approaches  that we know of, namely  those of  \cite{shen:2000} (his statistic is denoted by $\wt\alpha_n$)
 and \cite{KY2007}.  We discuss these approaches further in Section \ref{disc}, but  otherwise  restrict  discussion and analysis to $Q_n$.

The joint distribution of  the largest uncensored and the largest survival time in the sample is given in alternative ways in 
\cite{MR2020} and \cite{MRS2020}.
In the present paper we  apply the foundational  results in \cite{MRS2020} to obtain
exact finite sample as well as asymptotic  distributions  for 
 $Q_n$, 
which 
can be used  to assess whether  follow-up is sufficient in a sample.
The methods are illustrated with schematic and real data sets.

An important additional point is that  statistics such as 
$Q_n$ and $\wt\alpha_n$  can be used not only to  test for sufficient follow-up but 
 also to provide measures of how much follow-up there is in a sample.
Both these aspects are prominent in a paper by  \cite{Liu:etal:2018}   
where   testing for and measurement of sufficient follow-up in the TCGA pan-cancer clinical data resource are done on a very extensive scale 
in order to provide recommendations to cancer researchers wishing to assess the adequacy of clinical follow-up in a medical situation. 
\cite{Liu:etal:2018}
 processed  follow-up data files for 11,160 patients
across 33 cancer types, 
calculating median follow-up times as well as median times
to event (or censorship) based on the observed times for 
four endpoints (overall survival, disease-specific survival, disease-free interval, or progression-free interval).
They used  $Q_n$ and $\wt\alpha_n$  
to classify all $33 \times 4$ resulting  KMEs  as having sufficient or insufficient follow-up (or noted cases in which  tests were inconclusive)  in order to give  endpoint usage recommendations for each cancer type.
The analyses we present here help to validate the application of these tests in the data analysed in \cite{Liu:etal:2018}.

\section{Test Statistics for  Sufficient follow-up}\label{stes}
 
 \subsection{Notation and distributional setup}\label{ssD}
 For  the distributional results to follow we use the notation in  \cite{MRS2020},
which should be read in conjunction with the present paper.
We assume a general independent censoring model 
(``the iid censoring model") with right censoring.
A sample of size $n$ consists of observations on the sequence of
iid  (independent and identically distributed)  2-vectors 
$\big(T_i=T_i^*\wedge U_i, C_i={\bf 1}(T_i^*\le U_i);\, 1\le i\le n\big)$.
The $T_i^*$ with  continuous cumulative distribution function (cdf) $F^*$ on $[0,\infty)$ 
 represent the times of occurrence of an event under study, such as  the death of a person, the onset of a disease, the recurrence of a disease, the arrest of a person charged with a crime, the re-arrest of an individual released from prison, etc.
The $U_i$ with  continuous cdf $G$  on $[0,\infty)$  are censoring random variables, independent of the $T_i^*$.
In a sample of data from a population containing long-term survivors
we observe the  random variables $T_i=T_i^*\wedge U_i$, these being potential lifetimes  censored at a limit of follow-up represented for individual $i$ by the random variable $U_i$.
The   random variables
$ C_i={\bf 1}(T_i^*\le U_i)$ are censor indicators.
Let $M(n): =\max_{1\leq i\leq n}T_i $ be the largest observed survival time and let 
$M_u(n)$  be the largest observed {\it uncensored} survival time. 

The censoring distribution $G$ of the $U_i$ is always assumed proper (total mass 1), but  we allow the possibility that the distribution $F^*$ of the  $T_i^*$  is  improper.
We assume $F^*$  to be of the form
\begin{equation}    \label{FandF0}
    F^*(t)=pF(t),
\end{equation}
where $0<p\le 1$ and  $F$ is a proper distribution. We think of $F$ as being the distribution of susceptible  individuals in the population.
 Only susceptibles can experience the event of interest and have a potentially  uncensored failure time. 
The remainder of the population is immune to the event of interest or cured of it.
The presence of immunes is signalled by a value of $p<1$, in which case the  distribution $F^*$ is improper, with total mass $p$. 
Then $1-p$ is the proportion of immune or cured individuals in the population.
Observations on immunes are always censored; those on susceptibles may or may not be according as the corresponding $T_i^*>U_i$ or not.

Let $\Fbar^*(t)= 1-F^*(t)$, $t\ge 0$,
 denote the survival function (tail function) of $F^*$, and similarly for $\Fbar$ and  $\Gbar$.
Let $\tau_{F^*}= \inf\{t>0:F^*(t)=1\} $ 
(with the inf of the empty set equal to  $\infty$)  be the right extreme of the survival distribution $F^*$, 
and similarly $\tau_{F}$ and $\tau_{G}$
are  the right extremes of $F$ and $G$. 
 Let $H(t):=P(T_1\le t)$ be  the distribution of the observed survival times $T_i=T_i^*\wedge U_i$, with tail $\Hbar(t)=1-H(t)=P(T_i^*\wedge U_i>t   ) = \Fbar^*(t)\Gbar(t)$, $t\ge 0$, and right extreme $\tau_H=\tau_{F^*}\wedge \tau_G$. 
We always have $H(\tau_H)=1$, $G(\tau_G)=1$  and $F(\tau_{F})=1$. 
When $p=1$, so that $F^*\equiv F$, $F^*$ has total mass 1 and $\tau_{F^*}=\tau_{F}$;
when $p<1$ we have $\tau_{F^*}=\infty$, and  $\tau_{F}\le \tau_{F^*}$,
with the possibility that  $\tau_{F}< \tau_{F^*}$.

\subsection{Test statistics and procedure}\label{ssF}
 As test statistic  for sufficient follow-up we focus on  the statistic $Q_n$ proposed  in \cite{MZ1994}.
This is defined as follows.
Consider a sample of size $n$ with all survival times necessarily in $[0,M(n)]$, a number $N_u(n)$ of 
uncensored survival times, necessarily  in $[0,M_u(n)]$,
 a number  $N_c^<(n)$ of censored survival times  in $[0,M_u(n))$, and a number   $N_c^>(n)$ of censored survival times  in $(M_u(n), M(n)]$,
 thus   with  a total of $N_c(n)=N_c^<(n)+N_c^>(n)=n-N_u(n)$ censored survival times in the sample.
Set $\Delta_n:= 2M_u(n)-M(n)$.
As in  \cite{maller:zhoubook:1996}, p.81
we define
\be\label{qn2}
Q_n =\frac{1}{n} \#\{{\rm uncensored\ observations\ in}\ 
[\Delta_n, M_u(n)) \}.
\ee
(Note that we exclude $M_u(n)$ itself when counting the number of
uncensored  observations greater than $\Delta_n$.)
The statistic  $Q_n$ is   the proportion of uncensored observations in the interval $[2M_u(n)-M(n), M_u(n)]$, relative to the sample size $n$. It  measures  the length of the interval exceeding $M_u(n)$
but in a proportional rather than absolute way.
A rationale for the definition \eqref{qn2} is given in   \cite{maller:zhoubook:1996}, p.84.

The distribution of $Q_n$ was unavailable  when 
 \cite{maller:zhoubook:1996} was written and had to be simulated to get quantiles. 
 Our intention here is to get exact formulae 
(in Theorem \ref{thQ} below)
for  the distribution of $Q_n$ under the iid censoring model.
With these we can calculate asymptotic distributions (in Section \ref{dQ3})
and percentage points when estimates of $F$ and $G$ are made from data 
 (an example is in Section \ref{DA}).


Our test procedure will be as  follows. 
We have at hand survival  data with hypothesized cured individuals present 
and wish to test  for sufficient follow-up.
This is specified in  \cite{maller:zhoubook:1996}, p.81, to be the parametric condition  $\tau_{F}\le \tau_G$.
(For a rationale for this condition, see Sections 2.2 and 2.3 of   \cite{maller:zhoubook:1996}.)
We proceed by assuming the contrapositive hypothesis, $H_0: \tau_G<\tau_{F}$.
If $H_0$ is true the probability of seeing a large value of 
the test statistic  $Q_n$ is small.
So we will reject $H_0$ and conclude that   follow-up is sufficient if the observed value of the test statistic exceeds  a nominated quantile of its distribution under $H_0$.
 A test based on large values of $Q_n$ will reject  the hypothesis of insufficient follow-up with probability approaching 1 as sample size tends to infinity;
this follows from the asymptotic results in   Theorem \ref{thQ2}.
In order to understand the behaviour of $Q_n$, we first consider its finite sample properties.

\section{Understanding the sample properties of  $Q_n$}\label{dQ}


The value of $Q_n$ depends in a complicated way on the numbers of
censored and uncensored observations, the way they happen to occur below or above 
 $M_u(n)$,  and
 on the relative magnitudes of $M_u(n)$ and $M(n)$.
In order to calculate its distribution under the iid censoring model  we need to understand how it varies with these things. To do this we consider hypothetical sample situations, vary the mentioned quantities and see how the value of $Q_n$ changes. 

We begin by considering possible values of $\Delta_n$. 
We always have  $\Delta_n =2M_u(n)-M(n)\le  2M_u(n)-M_u(n)= M_u(n)$.
Possible values of $\Delta_n$ range from
$\Delta_n=- M(n)$ if $M_u(n)=0$, equivalently, if all observations are censored, to $\Delta_n=M_u(n)= M(n)$  if $M_u(n)=M(n)$,
 equivalently, if the largest observation is uncensored.
 We have $\Delta_n=0$ if it happens that $M_u(n)= M(n)/2$.
Thus we may have $\Delta_n<0$,  $\Delta_n=0$, or  $\Delta_n>0$.
When $\Delta_n\le 0$ then 
$[\Delta_n, M_u(n)) \supseteq [0,M_u(n))$ and
\eqref{qn2} gives $nQ_n=N_u(n)-1$.
At the other extreme, the interval $[\Delta_n,M_u(n))$ may be empty, and this is certainly so when  $\Delta_n=M_u(n)$.
 Whenever this occurs we set $Q_n=0$.
 
\begin{figure}
  \begin{subfigure}{7cm}
    \centering\includegraphics[width=6cm]{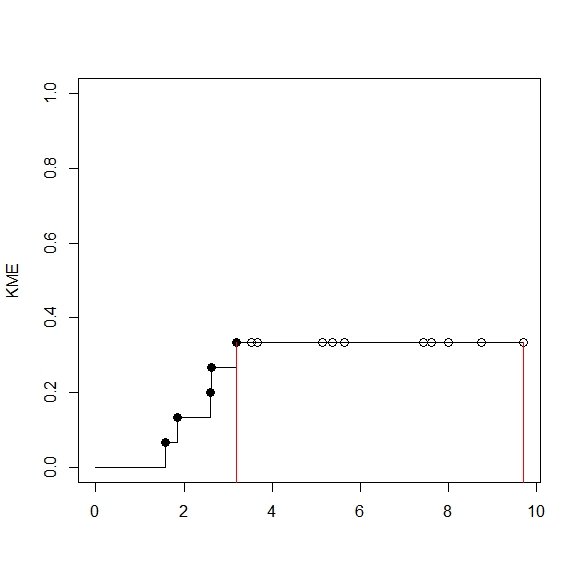}
    \caption{\it All censored observations $>M_u(n)$}\label{4figsa}
  \end{subfigure}
  \begin{subfigure}{7cm}
    \centering\includegraphics[width=6cm]{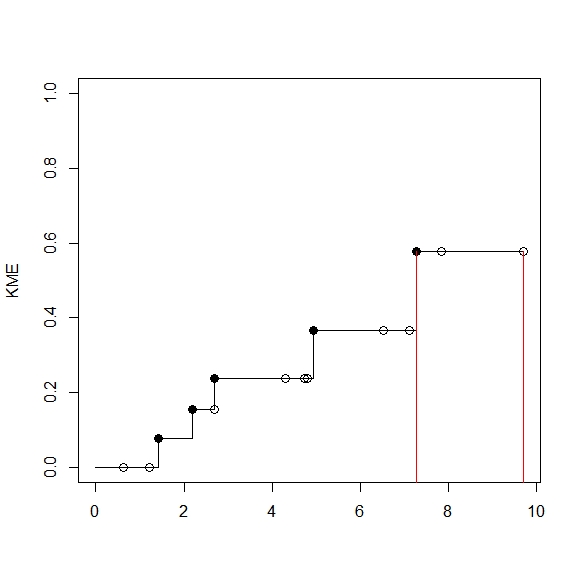}
    \caption{\it 2  censored observations $>M_u(n)$}\label{4figsb}
  \end{subfigure}
 
  \begin{subfigure}{7cm}
    \centering\includegraphics[width=6cm]{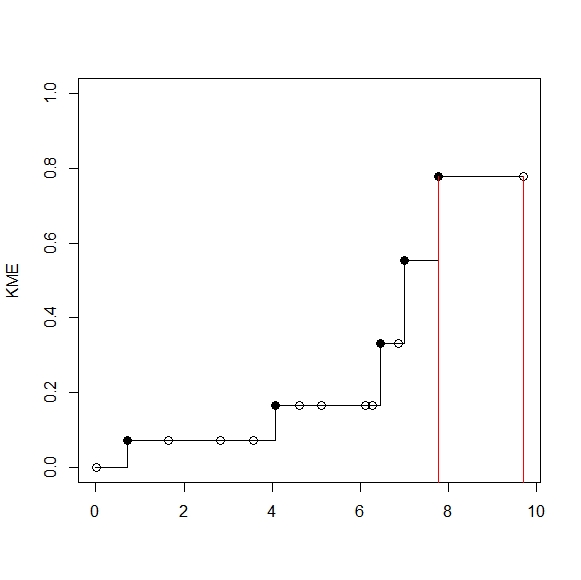}
    \caption{\it 1  censored observation $>M_u(n)$}\label{4figsc}
  \end{subfigure}
  \begin{subfigure}{7cm}
    \centering\includegraphics[width=6cm]{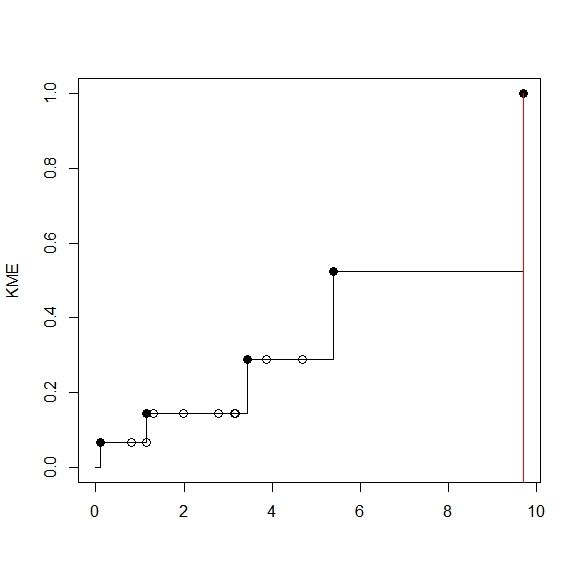}
    \caption{\it Largest observation uncensored}\label{4figsd}
  \end{subfigure}
  \caption{{\it  Schematic KME diagrams} }\label{4FIGS}
\end{figure}

 Now think of the way $Q_n$ changes if we 
 rearrange the conformation of the censored observations less than or greater than $M_u(n)$, by keeping $M(n)$ and $N_u(n)>0$ fixed and varying $N_c^<(n)$ and $N_c^>(n)$. 
It helps to visualise the various situations with schematic KME diagrams in the different cases, as we show in Figures \ref{4figsa} -- \ref{4figsd}.


 We start with an extreme case.
 \\
{\it Case 1:}\   $N_c^<(n)=0$, $N_c^>(n)>0$ (see Fig. \ref{4figsa}).
In this  conformation   all the censored observations in the sample form a  level stretch of the KME between 
 $M_u(n)$ and  $M(n)$. 
In this case  $M_u(n)$ takes the minimum possible value  for the sample under this kind of rearrangement,
$M(n)-M_u(n)$ takes the maximum possible  value,
$\Delta_n=2M_u(n)-M(n)= M_u(n) -(M(n)-M_u(n))$
takes the minimum possible  value,
and $Q_n$ takes the maximum possible  value  under this kind of rearrangement for the sample.
We reject $H_0:\tau_G<\tau_{F}$ and conclude there is sufficient follow-up if we observe large values of $Q_n$, so this arrangement accords with our intuition that a  (long)  level stretch on the KME between  $M_u(n)$ and  $M(n)$ indicates there is sufficient follow-up.
 \\
{\it Case 2:}\   $N_c^<(n)>0$, $N_c^>(n)>0$ (see Fig. \ref{4figsb}).
As censored observations are moved to the left of $M_u(n)$,
 $M_u(n)$ tends to increase and $M(n)-M_u(n)$  tends to decrease (it cannot increase).
 So  $\Delta_n$ will tend to increase and consequently
  $Q_n$ will tend to decrease. 
This accords with our intuition that a decrease in the  number of censored observations above
$M_u(n)$  and in the length of the level stretch of the KME between  $M_u(n)$ and  $M(n)$ makes it less likely to reject $H_0$, the hypothesis of insufficient follow-up. 
  
  Ultimately, continuing this process,  we reach: 
 \\
{\it Case 3:}\ $N_c^<(n)>0$, $N_c^>(n)=1$ (see Fig. \ref{4figsc}).
The one censored observation above $M_u(n)$ is $M(n)$ itself and 
$\Delta_n=2M_u(n)-M(n)$ will be close to or equal to $M_u(n)$.
The interval $[\Delta_n,M_u(n))$ is small and $Q_n$ is small, possibly equal to $0$ (this certainly occurs when  $\Delta_n=M_u(n)$).
 This accords with our intuition that a short  level stretch of the KME between  $M_u(n)$ and  $M(n)$ indicates via a small value of $Q_n$ that there is insufficient follow-up.

In these scenarios, $Q_n$ decreases monotonically from a sufficient follow-up situation to one with insufficient follow-up.

The actual values taken on by $Q_n$  in these scenarios depend  on the relative magnitudes of $M_u(n)$ and $M(n)$.
The possibilities are as follows.
Note that since $N_u(n)>0$, we have $M_u(n)>0$.

(a)\ When $0< M_u(n)\le \tfrac{1}{2} M(n)$, then $\Delta_n \le 0$ and  $[\Delta_n,M_u(n))\supseteq  [0,M_u(n))$.
 In this case 
 \ben
 Q_n =\frac{1}{n} \#\{{\rm uncensored\ observations\ other\ than}\ M_u(n)\}\ 
= \frac{N_u(n)-1}{n}.
\een
 This is the largest value $Q_n$ can take for a given sample.
 
(b)\  When $ \tfrac{1}{2} M(n)<  M_u(n)< M(n)$, then $\Delta_n>0$ and the interval
$[\Delta_n, M_u(n))$ contains, say, $k$ observations. 
We have $k\ge 0$ and $k\le n-1$ since there is at least one censored observation greater than $M_u(n)$, namely, $M(n)$.
So we can write
 \be\label{koc}
 Q_n =\frac{k}{n} =\frac{1}{n} \#\{{\rm uncensored\ observations\ in\ [\Delta_n, M_u(n))\} },
\ee
where $k$ decreases from its maximum value when $ M_u(n)$ is near $ \tfrac{1}{2}M(n)$, reaching 0 when
$ M_u(n)$ is near $ M(n)$.

There are also two other extreme cases to consider.

 (c)\  When $N_c^>(n)=0$, then $ M_u(n)=M(n)$, and the largest observation is  uncensored
 (see Fig. \ref{4figsd}). Then 
  $ \Delta_n= M_u(n)$, the interval 
  $ [\Delta_n, M_u(n))$ is empty,
 and  $Q_n=0$. 
 Here the level stretch has length 0 and  the low $Q_n$ value correctly reflects sufficient follow-up.
 (This case includes also the possibility that all observations are uncensored, corresponding to $N_u(n)=n$, and $k=n$.)
 But this Case (c)  means 
 there is no evidence of immunes and hence no issue of sufficient or insufficient follow-up. 
We condition on the non-occurrence of this event when calculating the distribution of $Q_n$.

 (d)\ When $N_u(n)=0$, all observations are censored, and, formally, $Q_n=0$. This anomalous or ambiguous case is of no interest and we condition on its non-occurrence also, 
   when calculating the distribution of $Q_n$.
%


\begin{figure}[h!]
\centering
\includegraphics[scale = 0.22]{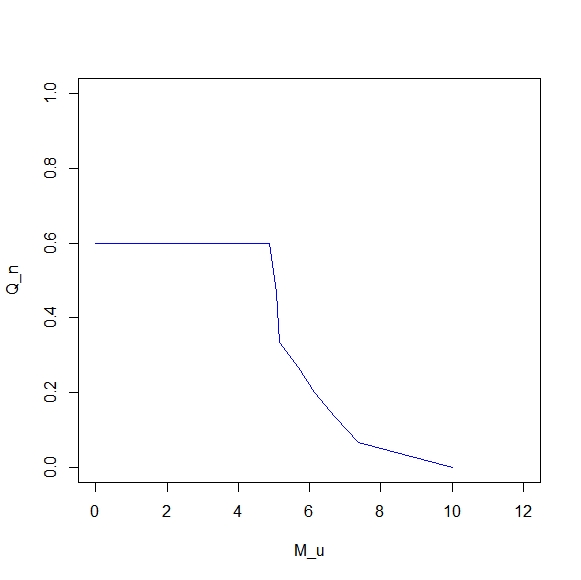}
\caption{{\it Possible Values for  $Q_n$} }\label{schem4}
\end{figure}

\section{Finite Sample  Distribution of  $Q_n$}\label{dQ2}
Given the discussion in the previous section,
in  calculating the distribution of $Q_n$
we will eliminate the cases $M_u(n)=0$ and
$M_u(n)=M(n)$, and partition the event of interest, $\{0<M_u(n)<  M(n)\}$, as
$\{0<M_u(n) \le  \tfrac{1}{2} M(n)\}
\cup \{ \tfrac{1}{2} M(n) < M_u(n)< M(n)\}$.
So we condition on $\{M_u(n)=t,M(n)=x\}$, where 
$0<t< x\le \tau_H$, and consider the  cases
$0<t\le \tfrac{1}{2} x$ and $\tfrac{1}{2} x < t <x$
separately.

Equivalently, we consider  {\bf Case A:}  $2t-x\le 0$,
 and   {\bf Case B:} $0<2t-x\le \tau_H$, with
 $0<t< x\le \tau_H$ in both cases. 
 Since we know the joint distribution of  $M(n)$ and $M_u(n) $  
 from Theorem 2.3 of \cite{MRS2020},
we can integrate to obtain 
the distribution of $Q_n$ conditional on 
$\{0<M_u(n)<  M(n)\}$.
One such calculation is carried out in the proof of  Theorem \ref{thQ}, and a large-sample version is in  Theorem \ref{thQ2}.

We set out some further preliminaries to Theorem \ref{thQ}.
Recall the formula \eqref{koc} for $Q_n$.
Throughout we keep $n>2$, $0<t< x\le \tau_H$ and $1\le r\le n-1$,  and will begin by 
conditioning on the event $\{M_u(n)=t,M(n)=x, N_c^>(M_u(n))=r\}$.
We need some separate notation in Cases A and B. 
 For Case A define
 \be\label{piA}
\pi^A(t):=
\frac{P(0<T_1^*\le t, T_1^*\le U_1) }{P(T_1^*\wedge U_1 \le t ) }
= \frac{\int_0^{t} \Gbar(y)\rmd F^*(y)}{H(t)}
\ee
(which does not depend on $x$)
and for Case B define
 \be\label{piB}
\pi^B(t,x):=\frac{P(2t-x< T_1^*\le t,T_1^*\le U_1) }{P(T_1^*\wedge U_1 \le t ) }
= \frac{\int_{2t-x}^{t} \Gbar(y)\rmd F^*(y)}{H(t)}.
\ee
Define also the probability
\be  \label{pctx}
p_c^>(t,x)=
 \frac{\int_{y=t}^{x} \Fbar^* (y)\rmd G(y)}
{ \int_{y=t}^{x} \Fbar^* (y)\rmd G(y)+H(t)},
\ee
and let
\be  \label{rhodef}
\rho^A(t,x):=(1-p_c^>(t,x))\pi^A(t)
\ {\rm and}\
\rho^B(t,x):=(1-p_c^>(t,x))\pi^B(t,x).
\ee

Define integers
$I_n^<:= \{i\in\N: T_i<M_u(n)\}$
and let
$\sigma^<$ be the smallest  $\sigma$-field 
making $(T_i,C_i)_{i\in I_n^<}$ measurable.
Likewise let $I_n^>:= \{i\in\N: T_i>M_u(n)\}$ and let 
$\sigma^>$ be the smallest  $\sigma$-field 
making $(T_i,C_i)_{i\in I_n^>}$ measurable.
A key result from Theorem 2.1 and Corollary 5 of 
\cite{MRS2020} is that, conditional on  the event
$\{M_u(n)=t, M(n)=x, N_c^>(M_u(n))=r\}$, or, equivalently,
 conditional on the event  $\{M_u(n)=t, N_c^>(M_u(n))=r\}$,
the   $\sigma$-fields $\sigma^<$ and $\sigma^>$ are
 independent, and the conditional probability of an event 
 $A^<$ in $\sigma^<$ can be calculated by substituting
truncated rvs $(T_i(t))$ 
having the distribution of $T_i$ given $T_i\le t$ 
for the $(T_i)$ in $A^<$.

With this setup we can now state  Lemma \ref{lem1}.
(Proofs of the lemma and the subsequent Theorems \ref{thQ} and  \ref{thQ2} 
are in the supplementary material.)

\begin{lemma}\label{lem1} {\bf Part (i):}\ We have for $1\le r\le n-1$,  $0<t<x\le \tau_H$, 
\be\label{L1}
P\big(N_c^>(M_u(n))=r\big|M_u(n)=t, M(n)=x\big)
=P\big(  Bin(n-2,\, p_c^>(t,x)) =r-1\big).
\ee
(With the lefthand side taken as  0 when $r=0$). 

{\bf Part (ii):}\
For $0\le k\le n-2$, 
\be\label{Qd}
P\big(nQ_n=k\big|M_u(n)=t,M(n)=x\big) =
P\big( Bin (n-2,\, \rho(t,x)) =k\big),
\ee
where $\rho(t,x)=\rho^A(t,x)$ in Case A 
and $\rho(t,x)=\rho^B(t,x)$ in Case B
(see \eqref{rhodef}).
\end{lemma}

We need one more formula:  by 
 Eq. (2.14) of \cite{MRS2020} 
we have
\bea\label{dtdx}
&&
P_n(\rmd t,\rmd x):= 
 P\big(M_u(n)\in \rmd t,M(n)\in \rmd x\big)\cr
 &&
 =
n(n-1) \Big(  \int_{y=t}^{x} \Fbar^*(y)\rmd G(y)+H(t)\Big)^{n-2} \,
 \Gbar (t)\rmd F^*(t)\,  \Fbar^*(x)\rmd G(x).\ \ 
 \eea
Next we can state 
Theorem \ref{thQ}.
 
\begin{theorem}\label{thQ}  
Assume the iid censoring model in Subsection \ref{ssD}.
Then for $n>2$, $k=0,1,2,\ldots, n-2$, 
\be\label{Qdis}
P\big(nQ_n=k\big| 0<M_u(n)<  M(n) \big) =
\frac{A_n(k)+B_n(k)}{D_n},
\ee
where
\be\label{Qdis3}
A_n(k)= 
 \int_{t=0}^{\tau_H/2} \int_{x=2t} ^{\tau_H}
 P\big( Bin (n-2,\, \rho^A(t,x))  =k\big)P_n(\rmd t,\rmd x)
\ee
and
\be\label{Qdis2}
B_n(k)=
\Big[\int_{t=0}^{\tau_H/2} \int_{x=t}^{2t} 
+ \int_{t=\tau_H/2}^{\tau_H} \int_{x=t}^{\tau_H} \Big]
P\big(  Bin (n-2,\, \rho^B(t,x))  =k\big)
P_n(\rmd t,\rmd x)
\ee
(recall \eqref{dtdx} for $P_n(\rmd t,\rmd x)$.)
The denominator in \eqref{Qdis} is
\bea\label{Qdis4}
D_n&=& P\big( 0<M_u(n)<  M(n) \big)\cr
 &=&
  1-
\left(\int_{t=0}^{\tau_H} \Fbar^*(z) \rmd G(z)\right)^n
-n\int_{t=0}^{\tau_H}H^{n-1}(t) \Gbar(t) \rmd F^*(t).
\eea
\end{theorem}

 \subsection{Probability mass functions of $Q_n$}\label{pmfs}
\begin{figure}[h!]
\centering
\includegraphics[scale = 0.5]{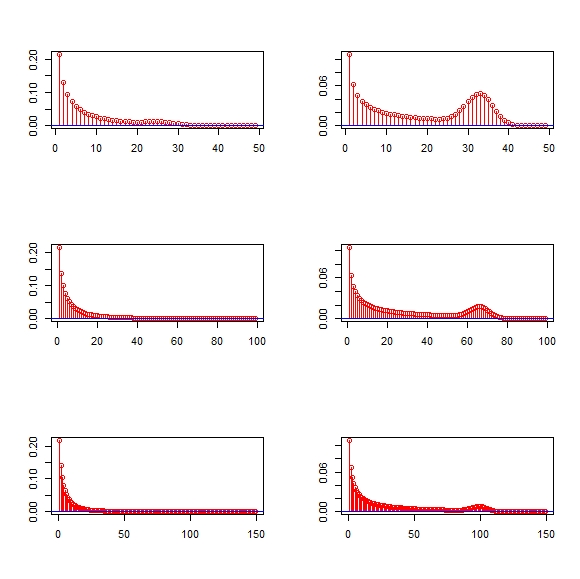}
\caption{{\it Probability mass functions for  $nQ_n$.
 $F=exp(1)$,  $p=0.8$.\\
Left column:  $G=U[0,3]$;  Right column: $G=U[0,6]$.\\
Top, middle, bottom row: $n=50, 100, 150$.
} }\label{fig2}
\end{figure}

 Figure \ref{fig2} 
 shows graphs of the probability mass functions (pmfs) 
 of $nQ_n$ calculated from \eqref{Qdis} for various scenarios  with $F$ exponential, $G$ uniform, and $n=50,100,150$.
For small $n$ the pmfs are bimodal, reflecting the two components on the RHS of  \eqref{Qdis}. The bimodality  is least prominent when censoring is heavy and disappears altogether as $n\to\infty$.
The pmfs in  Figure \ref{fig2}  are consistent with those obtained by  simulation in  Section 4.3 of  \cite{maller:zhoubook:1996}.
The tables in \cite{maller:zhoubook:1996} are based on 
 exponential survival distributions but remain relevant also for a 
certain scale family of distributions; see the Supplement for details.
 
%

We mentioned in connection with \cite{Liu:etal:2018} in Section \ref{intro} that $Q_n$ 
 can be used not only to  test for sufficient follow-up but 
also to provide a measure of how much follow-up there is in a sample.
In this respect the bimodality  evident in some  of the pmfs in  Figure \ref{fig2} is an unsatisfactory feature but as we pointed out it disappears  as $n\to\infty$ and samples of survival data are often of many thousands of individuals, as is the case in \cite{Liu:etal:2018}.

We go on to give the asymptotic distribution of  $Q_n$ in the next section.


\section{Asymptotic Distribution of  $Q_n$}\label{dQ3}
In this section we give the large sample
 distribution of  $Q_n$ in situations both of insufficient  (the main case of interest) and sufficient follow-up.
 We assume $G$ has a finite right endpoint $\tau_G$ and 
$\Gbar$ behaves linearly near $\tau_G$.
We also impose mild regularity conditions on $F$.
These conditions are satisfied when $G$ is uniform on $[0,\tau_G]$ and $F$ is exponential, for example. 
The asymptotic distribution of  $Q_n$ is shown to be geometric when $\tau_G<\tau_{F}$ and normal 
 when $\tau_{F}< \tau_G$, under these conditions.
Our main result for this section is:
 
\begin{theorem}\label{thQ2}[Asymptotic distribution of  $Q_n$]
Assume the iid censoring model and $0<p\le 1$ throughout.
We have the following limiting distributions in cases of interest. 
\newline
\noindent {\bf Case 1:}\   Assume 
 $\tau_G<\tau_{F}\le \infty$, so  $\tau_H=\tau_G<\infty$.  Suppose also 
\be\label{FGQ}
\Gbar(\tau_G-z)=a_G (1+o(1)) z \ {\rm as}\ z\dto 0,
\ee
for a constant $a_G>0$, 
and,   in addition, $F$  has a density  in a neighbourhood of $\tau_G$ which is positive and continuous at $\tau_G$. 
Then
\be\label{q00}
\lim_{n\to\infty}  P\big(nQ_n=k\big) 
= \frac{1}{4}\Big(\frac{3}{4}\Big)^{k},\ k=0,1,2,\ldots,
\ee
so $nQ_n$ is asymptotically geometric with parameter $1/4$. 

\noindent {\bf Case 2:}\  Assume  $\tau_{F}<\tau_G<2\tau_{F}$, so that  $\tau_H=\tau_{F}<\infty$.  
Suppose also that \eqref{FGQ} holds 
and in addition
\be\label{FGDspec}
  \Fbar(\tau_{F}-z)= a_{F}(1+o(1)) z, 
\ee
for a constant $a_{F}>0$. 
Then as $n\to\infty$ 
\be\label{q01}
\frac{ \sqrt{n}(Q_n-\nu^B)}{\sqrt{\nu^B(1-\nu^B)}}
\todr N(0,1),
\ee
where the parameter 
 $\nu^B= 
 p\int_{ 2\tau_{F}-\tau_G}^{\tau_{F}} \Gbar(y)\, \rmd F(y)/
(1-p\Gbar(\tau_F))
 \in (0,1)$.
\newline
\noindent {\bf Case 3:}\  Assume  $2\tau_{F}<\tau_G< \infty$, so that  $\tau_H=\tau_{F}<\infty$.  Suppose also that
 \eqref{FGQ} and \eqref{FGDspec} hold.
Then as $n\to\infty$ 
\be\label{q02}
\frac{\sqrt{n}(Q_n-\nu^A)}{\sqrt{\nu^A(1-\nu^A)}}
\todr N(0,1),
\ee
where  $\nu^A=p\int_{0}^{\tau_{F}} \Gbar(y)\, \rmd F(y)/
(1-p\Gbar(\tau_F))$.
\end {theorem}

\medskip\noindent{\bf Remarks.}\ (i)\
Case 1 with $\tau_G<\tau_F$ is a situation of insufficient follow-up, and in it $nQ_n$ has asymptotically a finite nondegenerate limit (a geometric rv).
Hence in this situation $Q_n\topr 0$ as $n\to\infty$, showing that the hypothesis of  insufficient follow-up will  be accepted in large samples (with probability approaching 1 as $n\to\infty$) when it is true.
When  follow-up is sufficient, i.e, Cases 2 and 3, $Q_n$ is ultimately normally distributed around positive levels $\nu^A$ or $\nu^B$ in large samples, and, depending on sample size, 
 the hypothesis of  insufficient follow-up will  be 
rejected, as it should be.
The specific formulae for the distributions in \eqref{q00}, \eqref{q01} and \eqref{q02}, enable the power calculations presented in the next section.

(ii)\
We remark that conditions \eqref{FGQ} and \eqref{FGDspec} are special cases of 
those imposed in Theorem 3.1 of \cite{MRS2020},
where more generally, a regularly varying function is
allowed in place of the linear factors in  \eqref{FGQ} and \eqref{FGDspec}.
However our present result is general enough for wide applicability.

\section{Power of the $Q_n$ Test}\label{Qpower}
In this section we use the asymptotic distributions of  $Q_n$ 
found under the assumptions  \eqref{FGQ} and the density condition on $F$ in Section \ref{dQ3},
 to calculate the power of the $Q_n$ test as the parameter $\tau_G$, reflecting the amount of follow-up, changes. 
In view of \eqref{q00}, it is more convenient to use $nQ_n$ than $Q_n$. 

We proceed by  calculating the 95-th quantile  $K_{0.95}$ of the  asymptotic distribution of  $nQ_n$ from \eqref{q00}, assuming  the hypothesis $H_0: \tau_G<\tau_{F}$ 
(insufficient follow-up) is true.
From \eqref{q00}
we can find $K_{0.95}$ explicitly as 
 \be\label{qan}
  K_{0.95} =   K_{0.95}(p,\tau_G) 
  =\frac{\log(0.05)}{\log (3/4)}-1= 9.41.
  \ee
Thus, under $H_0$, we have   $P(nQ_n>K_{0.95}) \approx 0.05$, for large $n$.
  Then we successively increase $\tau_G$ above $\tau_{F}$, 
hence in the region of the alternate hypothesis, and use
  \eqref{q01} and \eqref{q02} to calculate the corresponding values of $P(nQ_n>K_{0.95})$.
  Thus when  $\tau_{F}<\tau_G<2\tau_{F}$, according to   \eqref{q01}  we set
 \bea\label{q05}
 P(nQ_n>K_{0.95}) 
 &=&
 P\left(\frac{nQ_n-n\nu^B}{\sqrt{n\nu^B(1-\nu^B)}}
>  \frac{K_{0.95} -n\nu^B}{\sqrt{n\nu^B(1-\nu^B)}}\right)\cr
&\approx&
 P\left( N(0,1)>  \frac{K_{0.95} -n\nu^B} {\sqrt{n\nu^B(1-\nu^B)}}\right),
\eea
where 
 $\nu^B=p\int_{ 2\tau_{F}-\tau_G}^{\tau_{F}} \Gbar(y)\, \rmd F(y)/
(1-p\Gbar(\tau_F)))$;
and  when  $2\tau_{F}<\tau_G$,
 according to   \eqref{q02} 
 we replace $\nu^B$ in \eqref{q05}
by   $\nu^A= p\int_{0}^{\tau_{F}} \Gbar(y)\, \rmd F(y)/
(1-p\Gbar(\tau_F))$.
 
 Using $\nu$ to denote $\nu^A$ or $\nu^B$ as appropriate, we will use the approximation \eqref{q05}  only when $n\nu$ is large, and  take the  function of $\tau_G$ defined by 
 \be\label{Pdef}
 P(\tau_G;\nu): =  P\left( N(0,1)
> \min\left( \frac{K_{0.95} -n\nu}{\sqrt{n\nu(1-\nu)}}, 1.58\right)\right)
 \ee
 as an approximation to the power of the test.
 Keep $\tau_G>\tau_{F}$, and, at first, 
  $\tau_{F}<\tau_G<2\tau_{F}$.
 As  $\tau_G$ increases above $\tau_{F}$, 
 $\nu^B$ increases and $P(\tau_G; \nu^B)$ increases. 
 (Notice that $  K_{0.95}$ no longer depends on $\tau_G$ for values of $\tau_G>\tau_{F}$.)
 When $\nu^B=K_{0.95}/n$ then $P(\tau_G; \nu^B)$ reaches $0.50$, and once   $\tau_G$ reaches $2\tau_{F}$ then 
  $\nu^B=p\int_{0}^{\tau_{F}} \Gbar(y)\, \rmd F(y)/
(1-p\Gbar(\tau_F))$.
  For $\tau_G$ values greater than this $\nu^B$ is replaced in \eqref{Pdef} by $\nu^A= p\int_{0}^{\tau_{F}} \Gbar(y)\, \rmd F(y)/
(1-p\Gbar(\tau_F))$ and we note that 
  $\nu^A=\nu^B$ at the transition.
  For larger values of $\tau_G$,
  $P(\tau_G: \nu^A)$ stays constant at a value which  approaches 1 as $n\to\infty$. 
 
\begin{figure}[h!]
\centering
\includegraphics[scale = 0.15]{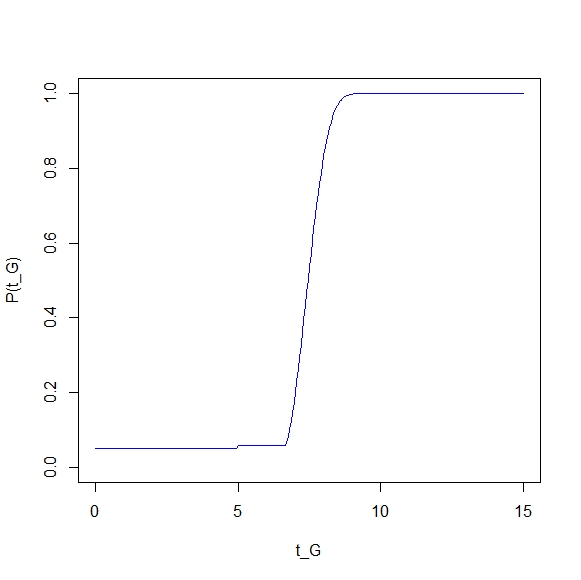}
\caption{{\it Power as a function of $\tau_G$ for  $Q_n$, with 
$F\sim exp(1)$, truncated at $\tau_{F}=5$; 
$G\sim [0,\tau_G]$; $n=100$.
} }\label{power}
\end{figure}

  In summary, the power function of the test appears to behave very well.
  We assume for illustration a sample size of $n=100$, for
$G$ a Uniform$[0,\tau_G]$ distribution, and for $F$ a unit exponential distribution truncated at a finite value $\tau_{F}=5$. Since  the probability in the tail of $F$ 
above 5 is less than $0.01$, this is effectively assuming 
 a unit exponential distribution for susceptible lifetimes.
  A graph of $P(\tau_G:\nu)$ for these parameter values is in Fig. \ref{power}.
  
\section{Dependent Censoring}\label{dep}
The distributional  results derived so far have been  based on the assumption of independence between the times of occurrence of the event under study and the censoring mechanism.
This assumption may not be tenable in some situations 
and there have been a number of studies where it has been relaxed.
See for example the competing dependent risks of leukaemia relapse and graft versus host disease analysed in \cite{kp} and \cite{kr}.
In order to assess the robustness of our results to departures from independence we consider the distribution of $Q_n$ in
a model where there is dependence between the survival and censoring distributions.

We assume  a sample consists of observations on the   2-vectors 
\ben
\big(T_i=T_i^*\wedge U_i,\ C_i={\bf 1}(T_i^*\le U_i);\, 1\le i\le n\big),
\een
where now the $T_i^*$ and  $C_i$ are dependent with a joint continuous distribution $H$
having marginal distributions $F^*$ and $G$ on $[0,\infty)$. 
We  model the dependence between $F^*$ and $G$  using a copula to
 connect the marginal distributions with the joint distribution. 
 
 Numerous copulas are defined and described  in \cite{Nel}, to which we refer for background.
For our bivariate setup we have 2 uniform random variables $W_1,W_2$, 
whose  joint distribution function is specified as
\be 
J(w_1,w_2,\theta):=P(W_1\le{w_1},W_2\le{w_2}),
\ee 
for a copula parameter  $\theta$ which quantifies the dependence between $W_1$ and $W_2$.
We restrict our discussion to  the class of {\it Archimedean copulas}  which contains 
many subfamilies  capable of representing different  dependency structures. 
The distribution function  of an Archimedean copula is written as:
\be 
J(w_1,w_2)=\Phi^{-1}\big(\Phi(w_1)+\Phi(w_2)\big),\  0\le w_1,  w_2\le 1, 
\ee
where the function  $\Phi$ is the {\it generator function } of the copula. 

We consider two generators which give rise to two Archimedean copulas: 
the  Frank and Ali-Mikhail-Haq (AMH)  copulas.
Each copula has an analytical expression that links its parameters to its related Kendall $\tau$ as a measure of association (\cite{smkr}).

The Frank copula (\cite{f}) has generator
\ben
 \Phi(t)=-\log\frac{e^{-\theta t}-1}{e^{-\theta}-1},\;\;\;\theta\in\R\setminus \{0\},
 \een
 giving rise to  the copula function
\ben
 J_{{\rm Frank}}(w_1,w_2)=-\frac{1}{\theta}\log\Big(1+\frac{(e^{-\theta w_1}-1)(e^{-\theta _2}-1)}
 {e^{-\theta}-1}\Big).
\een
The  Ali-Mikhail-Haq copula (\cite{amh}) has generator
\ben
 \Phi(t)=\log\Big(\frac{1-\theta(1-t)}{t}\Big),
 \;\;\;\; \theta\in[-1,1],
 \een
  giving rise to  the copula function
  \ben
   C_{{\rm AMH}}(w_1,w_2)=\frac{w_1w_2}{1-\theta(1-w_1)(1-w_2)}.
   \een

We proceed by assuming that  $H$ is a bivariate distribution with specified continuous marginals $F^*$ and $G$.
 By virtue of Sklar's theorem
 (\cite{sk}),  $H$ can be expressed in a unique way via a 2-copula $J$.
  In order to simulate an observation on 
  $(T_{i}^*, U_i)\sim H$, it is sufficient to simulate a vector $(W_1,W_2)\sim J$ 
  with values $w_1$ and $w_2$ 
  where the r.v.'s  $W_1$ and $W_2$ are uniform on $[0,1]$. 
Then
 \ben
 t^*=F^{*,\leftarrow}(w_1),\;\;u=G^{\leftarrow}(w_2),
 \een
is  an observation on  $(T^*,U)$ having the required joint distribution. 
(see \cite{smkr}, Appendix A).
 
\begin{figure}
\includegraphics[width=.2\textwidth]{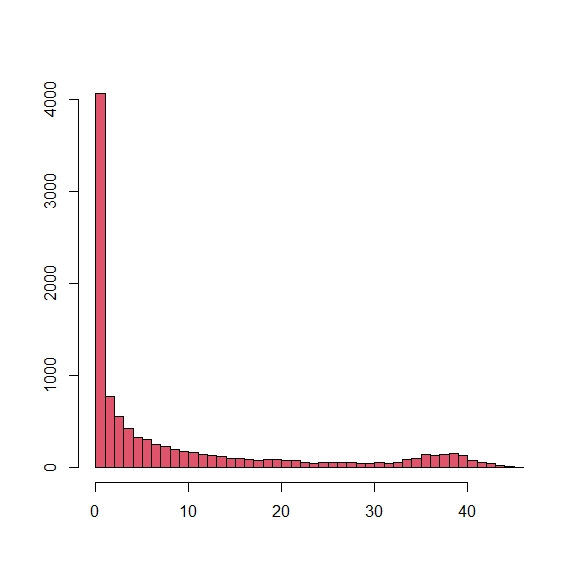}\hfill
    \includegraphics[width=.2\textwidth]{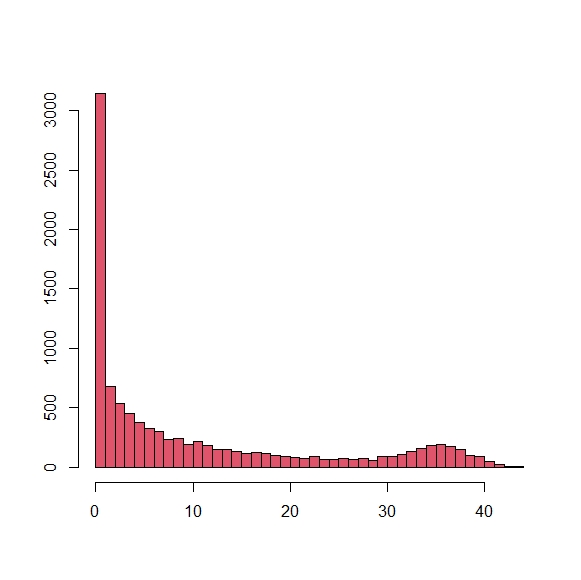}\hfill
    \includegraphics[width=.2\textwidth]{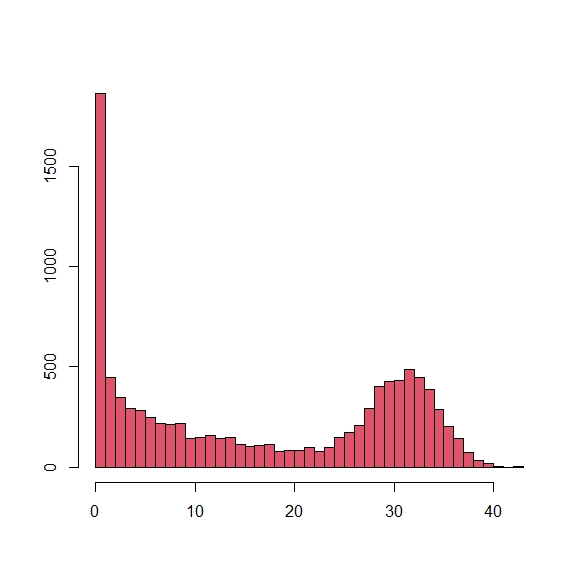}\hfill
     \includegraphics[width=.2\textwidth]{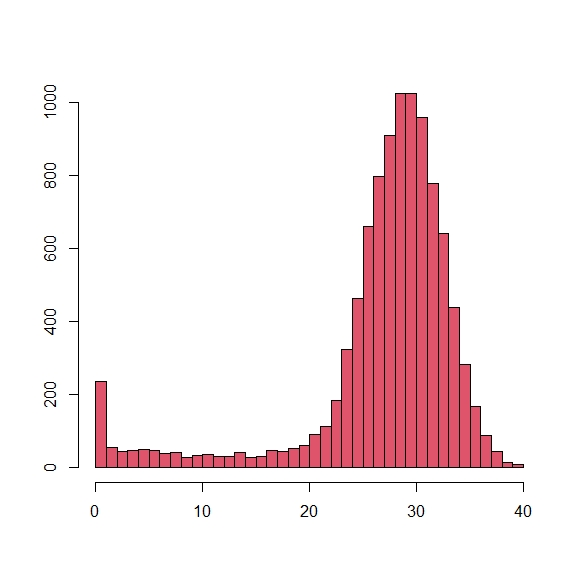}\hfill
    \includegraphics[width=.2\textwidth]{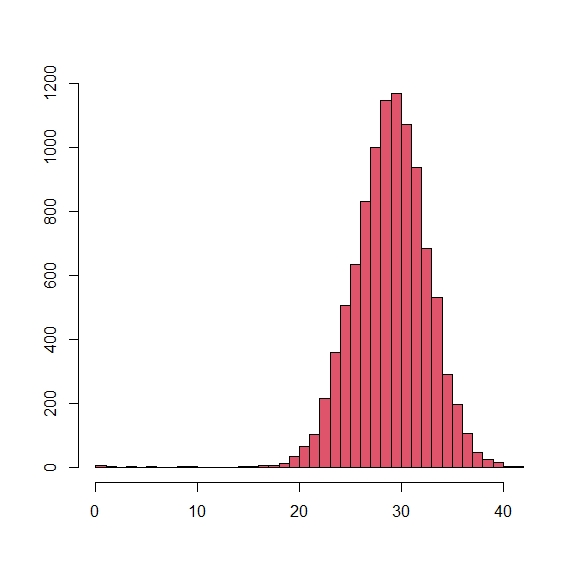}
    \\[\smallskipamount]
    \includegraphics[width=.2\textwidth]{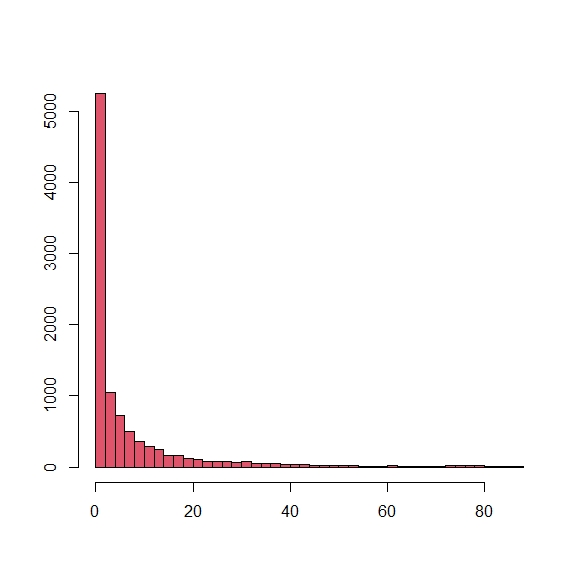}\hfill
    \includegraphics[width=.2\textwidth]{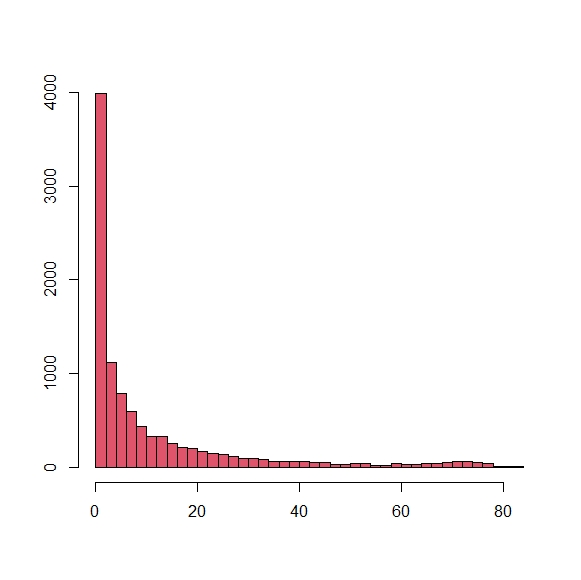}\hfill
    \includegraphics[width=.2\textwidth]{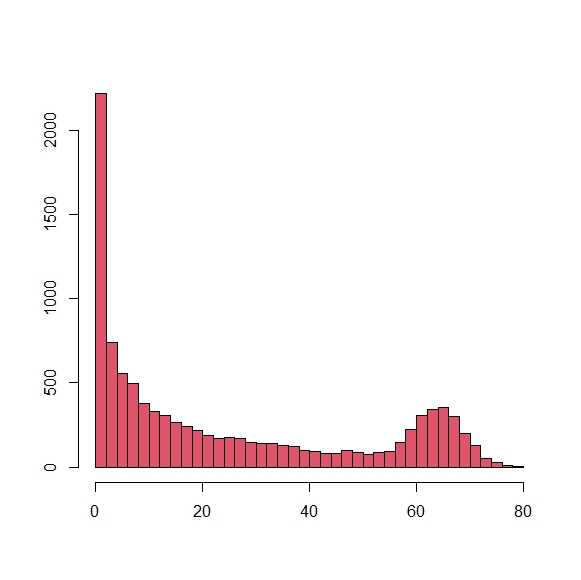}\hfill
    \includegraphics[width=.2\textwidth]{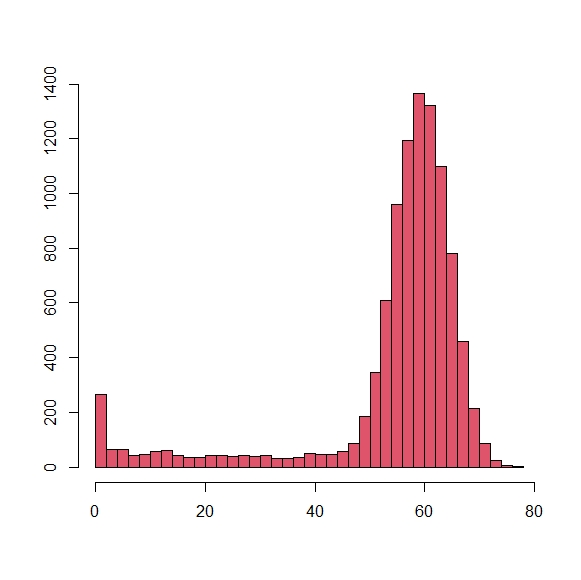}\hfill
    \includegraphics[width=.2\textwidth]{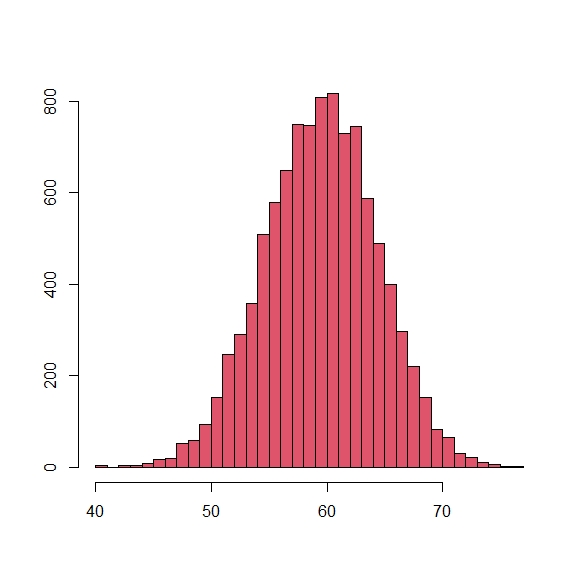}
    \\[\smallskipamount]
    \includegraphics[width=.2\textwidth]{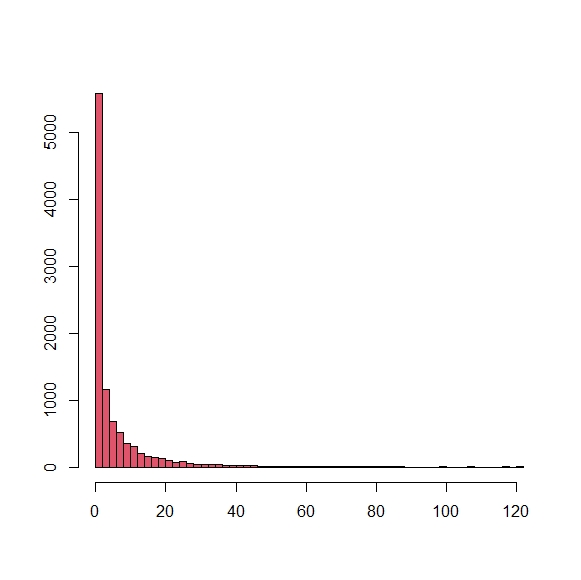}\hfill
    \includegraphics[width=.2\textwidth]{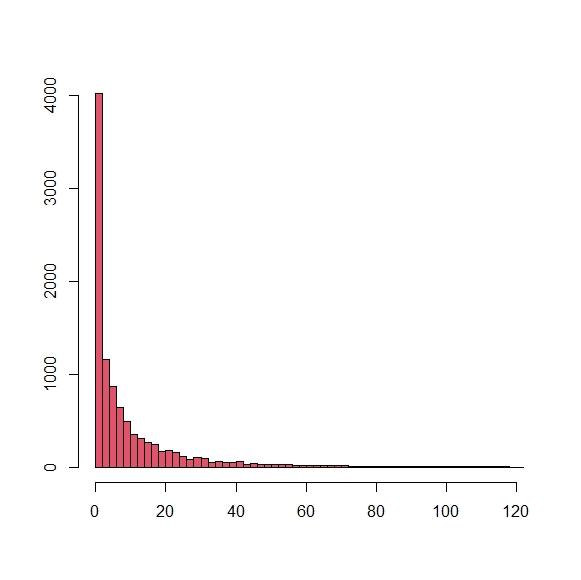}\hfill
      \includegraphics[width=.2\textwidth]{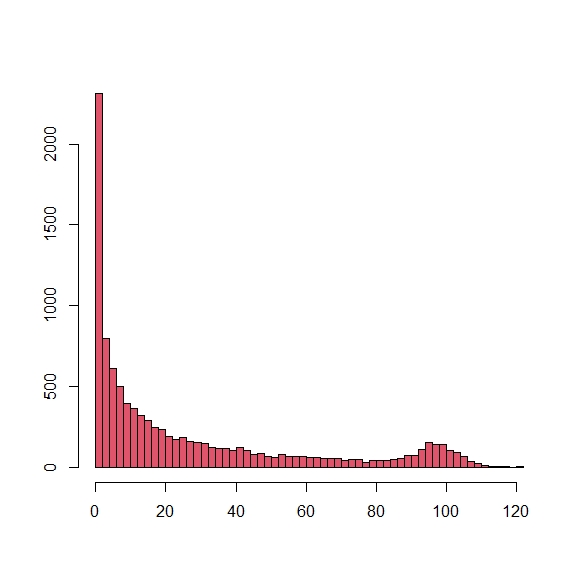}\hfill
      \includegraphics[width=.2\textwidth]{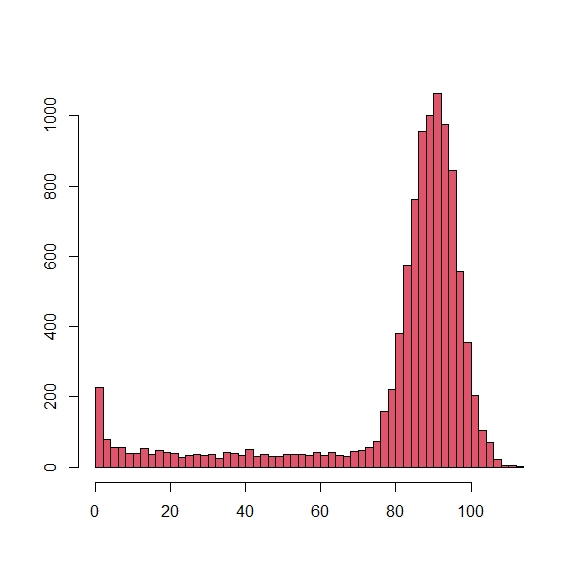}\hfill
    \includegraphics[width=.2\textwidth]{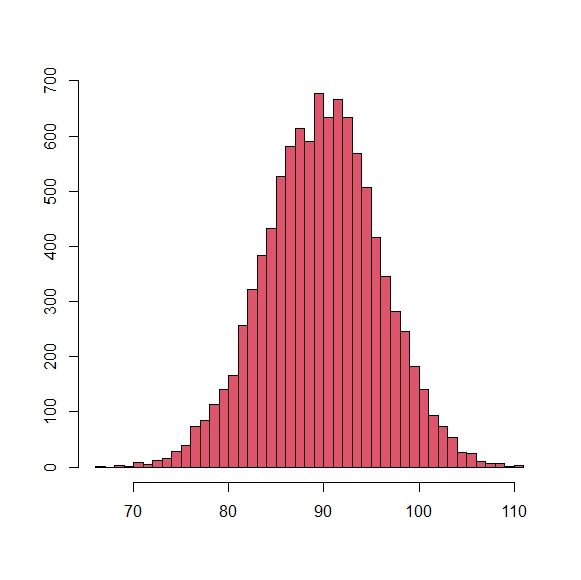}
      \caption{
      Probability mass functions for  $nQ_n$ with  Frank copula for dependence.
 $F=exp(1)$,  $p=0.8$, $G=U[0,6]$.
Top, middle, bottom panel: $n=50, 100, 150$.
Left to right:  $\theta=300,6,0,-6,-300$. 
      }\label{fig:foobar}
\end{figure}

\begin{figure}
    \includegraphics[width=.2\textwidth]{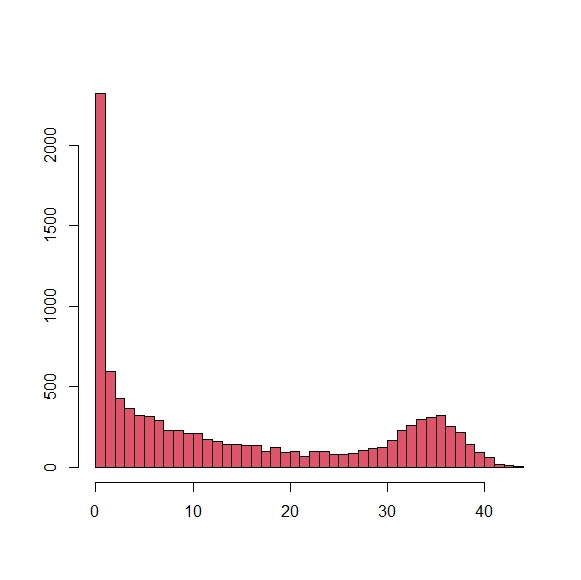}\hfill
    \includegraphics[width=.2\textwidth]{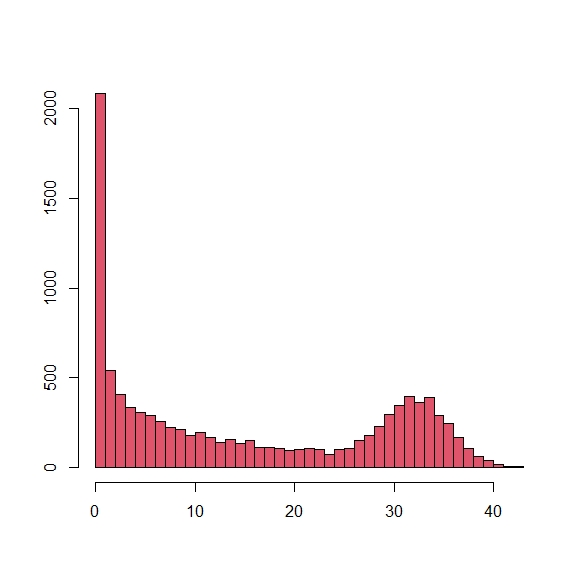}\hfill
    \includegraphics[width=.2\textwidth]{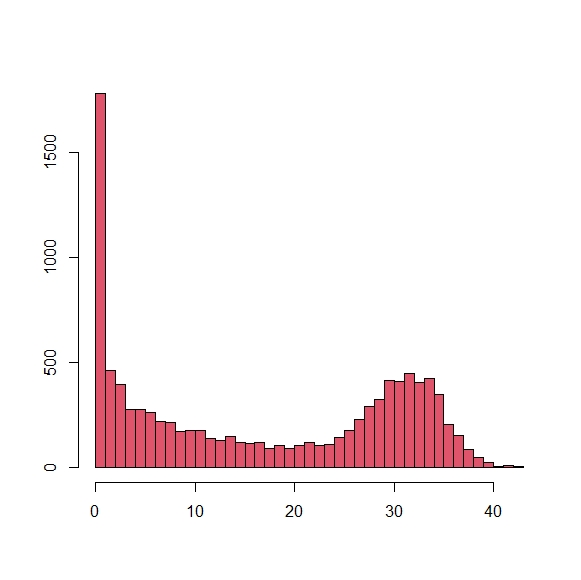}\hfill
     \includegraphics[width=.2\textwidth]{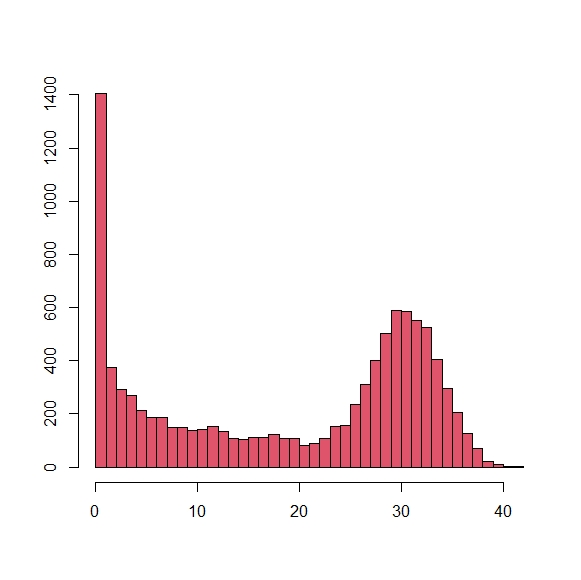}\hfill
    \includegraphics[width=.2\textwidth]{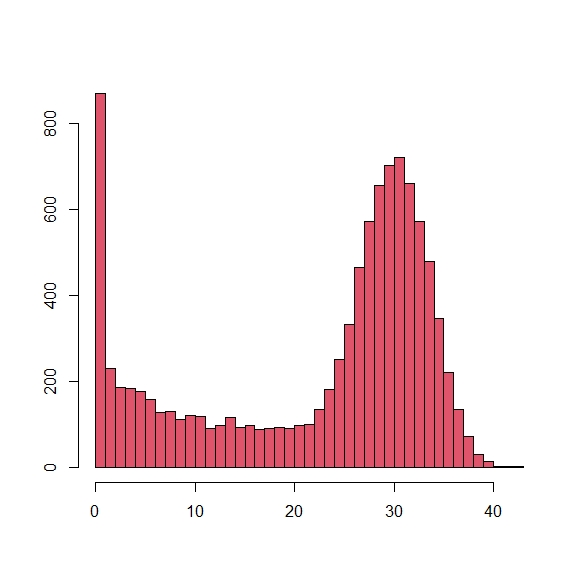}
    \\[\smallskipamount]
    \includegraphics[width=.2\textwidth]{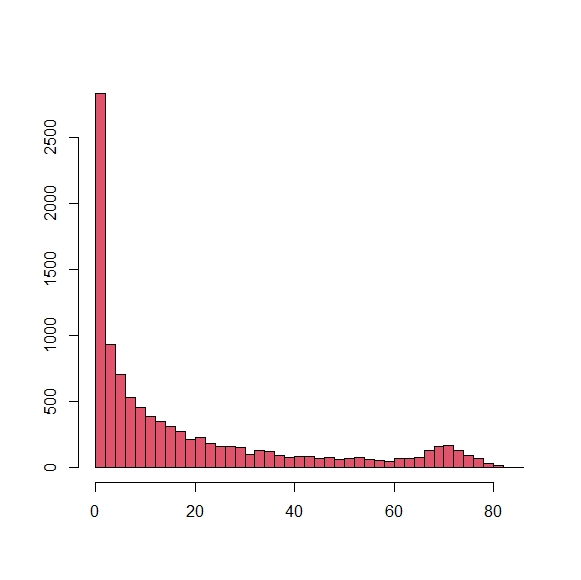}\hfill
    \includegraphics[width=.2\textwidth]{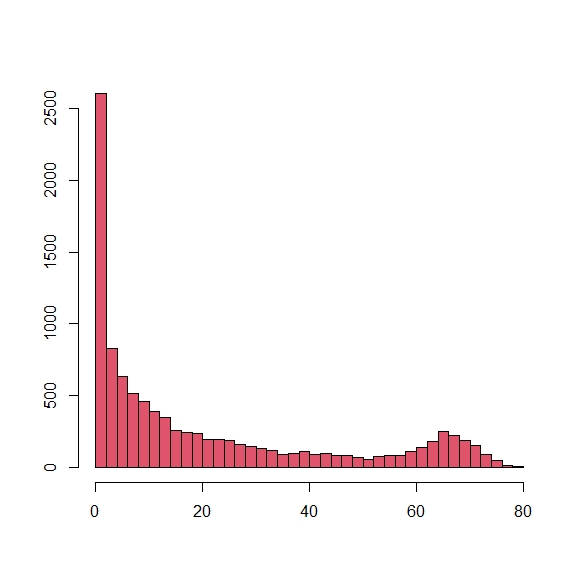}\hfill
    \includegraphics[width=.2\textwidth]{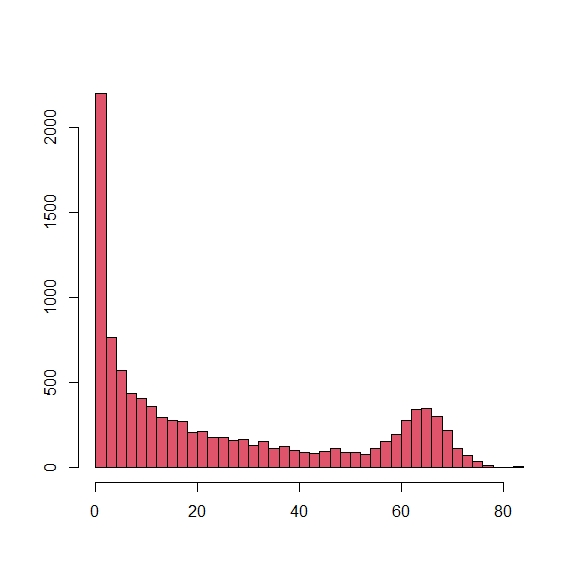}\hfill
    \includegraphics[width=.2\textwidth]{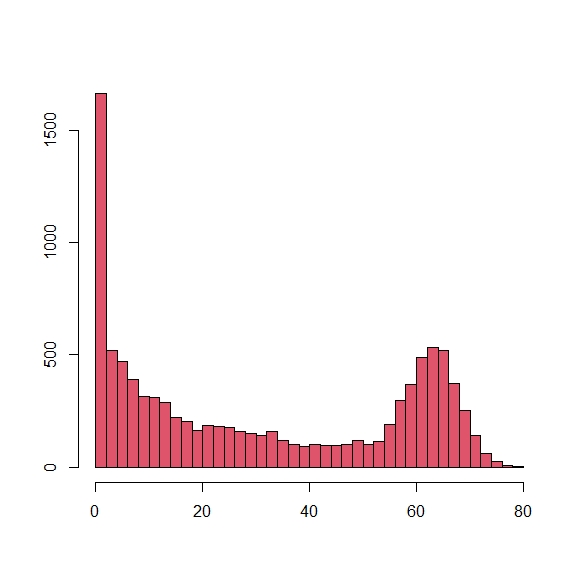}\hfill
    \includegraphics[width=.2\textwidth]{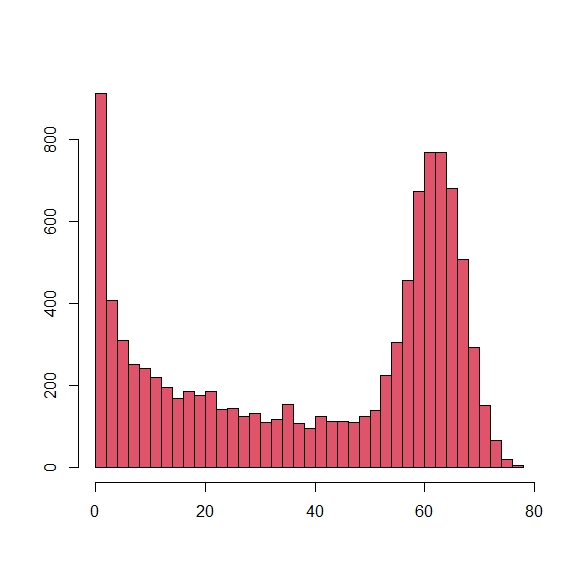}
    \\[\smallskipamount]
    \includegraphics[width=.2\textwidth]{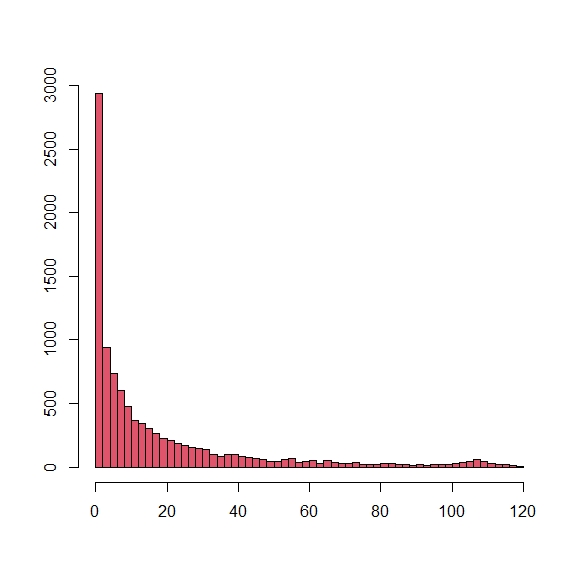}\hfill
    \includegraphics[width=.2\textwidth]{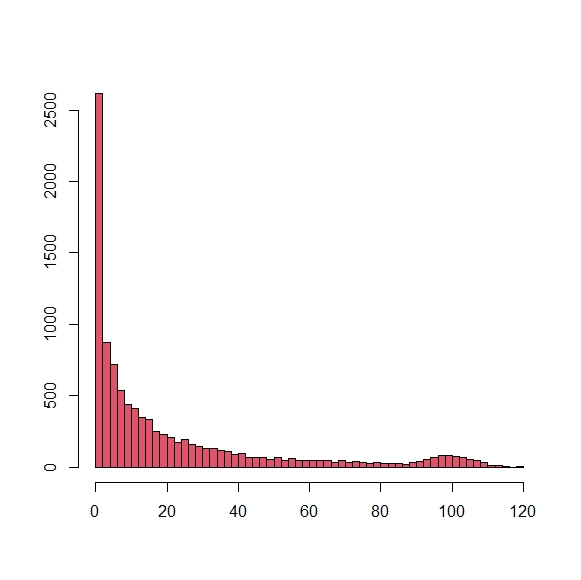}\hfill
      \includegraphics[width=.2\textwidth]{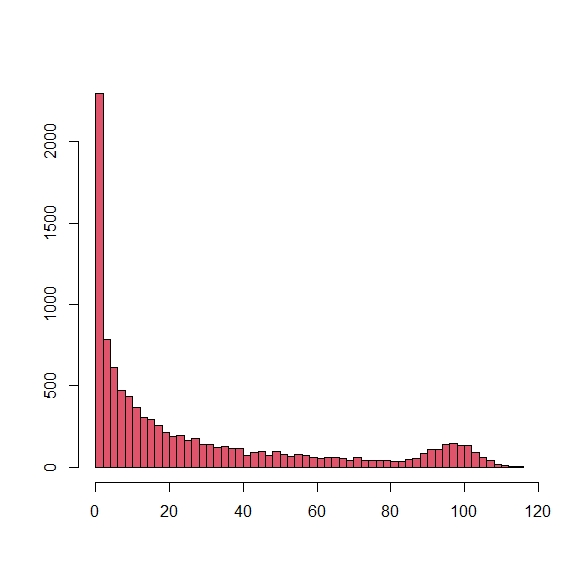}\hfill
      \includegraphics[width=.2\textwidth]{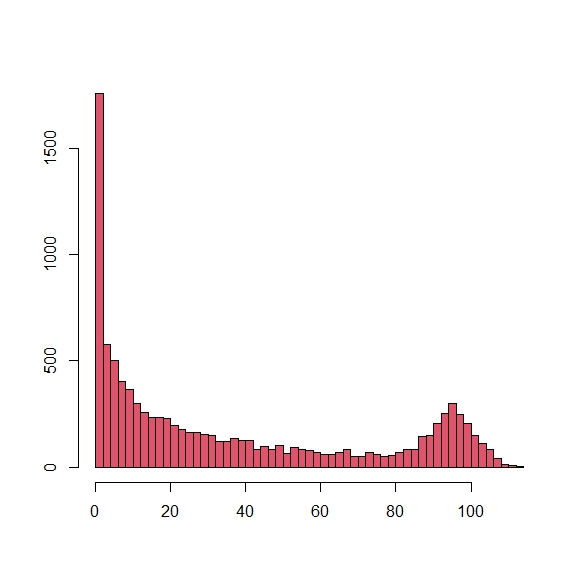}\hfill
    \includegraphics[width=.2\textwidth]{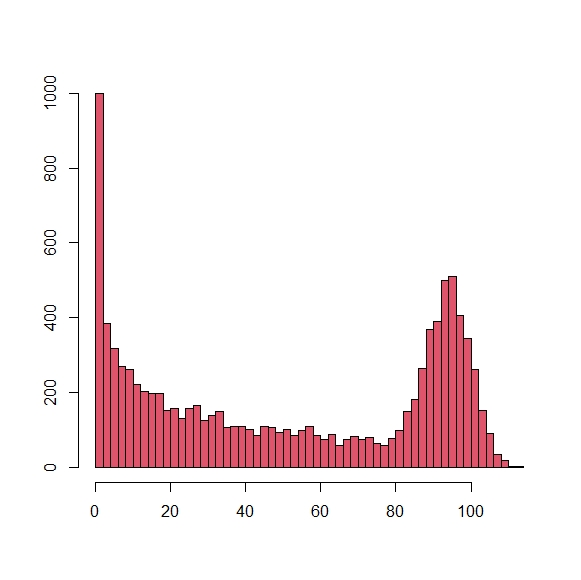}
  
      \caption{
            Probability mass functions for  $nQ_n$ with AMH copula for dependence.
 $F=exp(1)$,  $p=0.8$, $G=U[0,6]$.
Top, middle, bottom panel: $n=50, 100, 150$.
Left to right:       $\theta=1, 0.5, 0, -0.5, -1$.
      }\label{fig:amhbar}
\end{figure}

 For our robustness analysis  we simulated samples of size 
   $n=50,100,150$,  from a $J$ based on the Frank and AMH copulas
 for various values of $\theta$.
We took  $F^*=pF$, where $F$ is  exponential  with parameter 1, $G=U[0,6]$ and $p=0.8$.
In each sample  we calculated the value of  $Q_n$ and  repeated
  this $N=10000$ times to draw up the pmfs of  $Q_n$
  (Figures \ref{fig:foobar} and \ref{fig:amhbar}).
  In each figure the pmf for $\theta=0$ corresponds to independence. 
  Viewing from  the centre panel left we see that introducing negative dependence tends to concentrate the mass near small values of $Q_n$;
    viewing from  the centre panel right shows that posiative dependence tends shifts the pmfs closer to normal. Percentage points calculated from these distributions could be used for correction if dependence is assumed or detected in a sample.

\section{Data Example}\label{DA}
\begin{figure}[h]\label{KMEsfit}
  \centering
\hskip -1cm
  \begin{subfigure}{5cm}
\includegraphics[width=6cm, height=6cm]{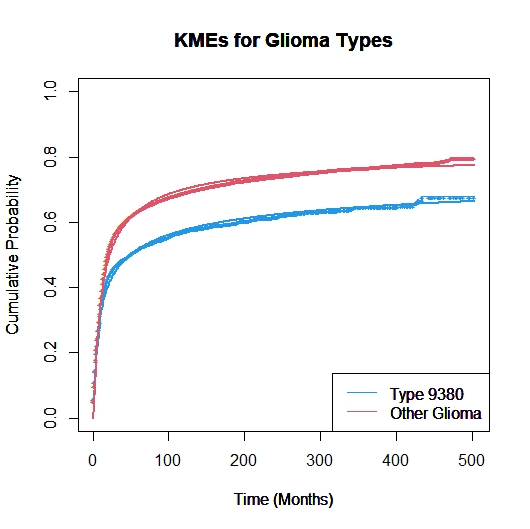}
  \caption{\it  Survival KMEs  for Glioma Types with fitted generalised gamma  distributions.}\label{uncensKMEs}
  \end{subfigure}
\hskip2cm 
  \begin{subfigure}{5cm}
    \centering\includegraphics[width=6cm, height=6cm]{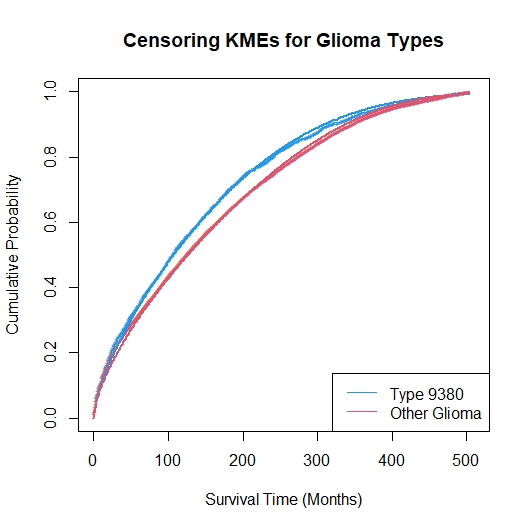}
\caption{\it Censor KMEs for  Glioma Types  with fitted generalised gamma  distributions.}\label{censKMEs}
  \end{subfigure}
\end{figure}

We calculated    
KMEs for the survival and censoring distributions from  data for glioma (brain cancer) patients contained in  the  \cite{SEER:2019} database   (Figures \ref{uncensKMEs} and  \ref{censKMEs}).
This is observational rather than from a randomised clinical trial but we use it here just to illustrate how consideration of followup with the $Q_n$ test can add value to an analysis of censored survival data. 

The data was subdivided into two classes: Type 9380 (4248 patient records, of which 2075 are censored) and Other Types\footnote{Types 9381, 9382, 9383, 9384, 9392, 9401, 9430, 9451, 9440, 9441, 9442, as used in
\cite{YXGJHC}.} (54375 patient records, 20482  censored) according to the SEER classification scheme.

%


%


The maximum values of the KMEs are 0.68 for Type 9380 and 0.80 for Other Types,
 suggesting cured probabilities of 0.32 for Type 9380  and  0.20 for Other Types.
These are significantly different from 0 based on Greenwood's formula for the variance of the KME, or using the tables in \cite{maller:zhoubook:1996}.
But is followup sufficient for us to be confident in concluding cured components?
For Type 9380 we note a level stretch at the end of the KME at 503 months, with the largest uncensored observation at 435 months. 
Measuring back $503-435= 68$ months we come to 367 months, and the number of uncensored observations in $[367,435)$ is found to be 3.   So $nQ_n=3$ for this data.  Since the sample size is large we can use
Case  1 of Theorem \ref{thQ2} to calculate the p-value of the statistic for testing
 $H_0: \tau_G<\tau_F$ as
$(3/4)^3= 0.42$ which is far from significance.
So we do not reject $H_0$  and decide that  followup is insufficient for this data.
 The same conclusion follows using the tables in \cite{maller:zhoubook:1996}.
We can see that, despite apparent levellings in Figure \ref{uncensKMEs}, the curves continue to rise, with late deaths still occurring up till 400 months after diagnosis. 

This test is nonparametric. But having failed to  reject $H_0: \tau_G<\tau_F$,  Cases 2 and 3 of Theorem \ref{thQ2} become relevant, and we can also consider the possibility that   $\tau_G=\tau_F=\infty$. The latter case is not included in Theorem \ref{thQ2}   as the corresponding result requires assumptions on the asymptotic behaviour of $(M_u(n), M(n))$, 
  but it can be shown that $Q_n$ is asymptotically normal under reasonable assumptions.
In order  to apply  Theorem \ref{thQ2}, estimates of the distributions of $F$ and $G$ are needed. 
Generalised gamma distributions with a cured component 
as discussed in \cite{Jackson:2016} and \cite{Amdahl:2020}\footnote{For alternative nonparametric fitting see \cite{RJ-2021-027}.} 
fit the data quite well  (Figures \ref{uncensKMEs} and  \ref{censKMEs}) and can be used for further 
analysis. Knowledge of the exact distribution allows for simulation via MCMC, for example.

For  the Other glioma types the value of $nQ_n$ is 13 and the p-value using Case  1 of Theorem \ref{thQ2}
is 0.02. For this data we reject $H_0: \tau_G<\tau_F$  and decide that  followup is sufficient.
There is only a small level stretch at the end of the KME, but the conformation of  censored and  uncensored observations mean that  it is significant.
Further details of the data and analysis are in the Supplement.

This brief discussion is not meant to be a substantial analysis of this data, at all, but it highlights the valuable information present in the observations at the right hand end of a KME, and how the $Q_n$ statistic can guide us in interpreting it.

\section{Discussion}\label{disc}
Calculation of exact distributions (under the iid censoring model) makes unnecessary the need for simulations of percentage points, though 
in  practice the unknown distributions 
must be estimated or postulated.
This includes making assumptions about their right hand endpoints, and,  especially, whether they are finite or not. 

Most practical is to assume $\tau_G<\infty$ since observation must always cease at some finite point.
In many cases the assumption $\tau_{F}<\infty$  may also be natural.
Certainly in real survival data no individual lives forever, but we would have $\tau_{F}=\infty$  for example when studying the occurrence of an infectious disease where an immune individual would never contract the disease no matter how long the follow-up. 
(For just such an analysis with children immune to malaria, see 
 \cite{malaria}.)
Regardless of the situation,  in modelling exercises it is not uncommon to use  an exponential,  Weibull, lognormal,  or Gumbel,
with infinite right endpoints,  as the lifetime distribution.
In doing so we accept that the probability of seeing an extremely long lifetime under the assumed model is negligible, so the theoretical  approximation is good enough for practical purposes. 

It is natural to base tests  for sufficient follow-up on the length of the  interval
$(M_u(n), M(n)]$,  or the  number of censored survival times larger than the largest uncensored survival time, or some combination or variant of these.
Kaplan-Meier plots provide very strong intuition in this respect; see for example the very evocative plots in \cite{powlesetal.2021}, 
who are expressly concerned with 
``determin(ing) which patients $\ldots$ are cured after surgery."

As other test statistics we might use the difference between the extremes, $M(n)-M_u(n)$,  or a standardised version of this such as  $R_n=1- M_u(n)/M(n)$, which is in $(0,1)$. Formulae for  their distributions are in \cite{MRS2020}. 
These variables measure the absolute or relative length of the level stretch of the KME rather than a proportion of observations related to them, as $Q_n$ does. 
At present their properties remain to be investigated  in detail. 
We note that they, like $Q_n$, are very sensitive to the occurrence of one or a few failures in the righthand end of the KME. This is a robustness issue such as has to be addressed in any statistical analysis. 
A test for outliers in the iid model  is in \cite{MZ1994}.
The  \cite{shen:2000} statistic $\wt\alpha_n$
is suggested by similar arguments to the way $Q_n$ was obtained in  \cite{maller:zhoubook:1996}, though the rationale seems not so obvious and the sample properties not so clear as for $Q_n$. 
 \cite{shen:2000} does not give a formula for  the  distribution of $\wt\alpha_n$. We expect that one could be obtained as a modification of the 
way we calculated the distribution of  $Q_n$ in  Theorem \ref{thQ}, but the computations are more difficult. 
Shen reports, based on some limited  simulations, that $\wt\alpha_n$  sometimes performed better than $Q_n$ in terms of Type 1 error and power.
A more extensive investigation of this is warranted.
 
 \cite{KY2007} give a succinct 
 overview of the $Q_n$ statistic as presented in \cite{maller:zhoubook:1996} and discuss some of its properties.
 They draw attention to a perceived deficiency of the statistic,
 as follows. Suppose follow-up in a sample were hypothetically extended beyond what's currently there.
It's possible then that a  long-lived susceptible individual
with  currently censored lifetime 
 may die during the extended follow-up period, and that the $Q_n$ value calculated on the extended sample then decreases from its former value, possibly even to 0.
 Klebanov and Yakovlev see this non-monotone behavior as problematic.  But there is really no contradiction here. 
On the extended sample, with its late failure, the new, low, $Q_n$ is correctly registering that there is insufficient follow-up.  If in this hypothetical situation we continue to increase follow-up, susceptible individuals will continue to die (all do, eventually), while those who are immune will remain so, and with extended follow-up. 
The $Q_n$ value may well  fluctuate, but eventually only cured individuals will be left and the continued follow-up will give rise to increasing values of $Q_n$,  till it reaches its maximum value for the sample. At each stage $Q_n$ is correctly (according to its constitution) indicating the extent of follow-up.  
There is no reason why $Q_n$ should be monotone in a  hypothetical situation of increasing follow-up. 

Nevertheless the quite different approach  in  \cite{KY2007} provides a potentially useful  perspective on the problem.
Unfortunately however their proposed statistic is technically complex and very non-intuitive and its application would likely be restricted to specialist statisticians, whereas an approach based on the length of the level stretch at the end of the KME is  highly visible and interpretable, and easily understood by practitioners.  Shen's statistic is also constructed in this way.
That there is a need and a desire for a summary statistic with these properties  is well exemplified by the \cite{Liu:etal:2018} analysis.

Understanding how $Q_n$ depends on sample properties of censored data, and the formulae for the exact and  asymptotic distributions of  $Q_n$ we have obtained, open the way to its more general use in the analysis of survival data with immune or cured individuals.
We note that under  $H_0: \tau_G<\tau_F$, the hypothesis of insufficient follow-up, with some reasonable side conditions, the asymptotic distribution of $Q_n$ is completely non-parametric (cf. \eqref{q00}). 
Future directions of research could usefully include issues of sufficient follow-up in competing risks analysis, and in multivariate survival analysis with cured individuals. 
For the latter, see \cite{CS1}, \cite{CS2} and \cite{CBAM}.

  \cite{shen:2000} and 
 \cite{KY2007} 
 quote  \cite{maller:zhoubook:1996} to the effect that $Q_n$ is by no means the last word on the subject, and this is worth stressing again here. 
Having formulae for the exact and asymptotic distributions of $Q_n$ in the iid censoring model is a big step forward, but there is much still to be explored in the analysis of sufficient follow-up and cure models in general. 

\medskip\noindent{\bf Acknowledgments}\ 
We are pleased to thank  Muzhi Zhao for helpful comments.

\begin{onehalfspacing}

\bibliography{bibross.bib}

\end{onehalfspacing}

\end{document}